\documentclass[12pt]{amsart}
\usepackage{amscd, amssymb}
 
\newcommand{\bdes}{\begin{description}}
\newcommand{\edes}{\end{description}}
\newcommand{\beqn}{\begin{eqnarray*}}
\newcommand{\eeqn}{\end{eqnarray*}}
\newcommand{\PP }{{\mathbb P}}
\newcommand{\QQ }{{\mathbb Q}}
\newcommand{\CC }{{\mathbb C}}
\newcommand{\ZZ }{{\mathbb Z}}

\newcommand{\NN }{{\mathbb N}}
\newcommand{\ke }{{\varepsilon }}

\newcommand{\ka }{{\alpha}}

\newcommand{\kd }{{\delta}}
\newcommand{\kl }{{\lambda}}

\newcommand{\hb }{{\hbar}}

\newcommand{\hbp}{\hbar\frac{\partial}{\partial t_i}}
\newcommand{\kline}{[\text{line}]}                          

\newcommand{\Lie }{{\mathrm{Lie}}}
\newcommand{\ad }{{\mathrm{ad}}}

\newcommand{\ti }{\times}
\newcommand{\ot }{\otimes}
\newcommand{\ra }{\rightarrow}

\newtheorem{theorem}{Theorem}

\newtheorem{def/th}{Definition/Theorem}

\newtheorem{lemma}{Lemma}

\begin{document}
\author{Bumsig Kim
\\ University of California at Davis}
\title{Quantum Hyperplane Section Theorem For Homogeneous Spaces}
\date{June 3, 1997, Revised December 31, 1997}
\begin{abstract}
We formulated a mirror-free approach to the mirror conjecture, namely,
quantum hyperplane section conjecture, and proved it in the case of
nonnegative complete intersections in homogeneous manifolds.
For the proof we followed the scheme of Givental's proof of a
mirror theorem  for toric complete intersections.
\end{abstract}
\maketitle
\pagestyle{plain}

\section{Introduction}

Quantum cohomology of a symplectic manifold 
is a certain deformed ring of the ordinary
cohomology ring with parameter space given by the second cohomology
group.  It encodes enumerative geometry
of rational curves on the manifold. In general it is difficult to compute the
quantum cohomology structure. On the other hand,
mirror symmetry predicts an answer to
traditional questions of counting the virtual numbers of
rational curves of a given degree on a three-dimensional Calabi-Yau
manifold, which amounts to knowing the
structure of the quantum cohomology.
In the large class of Calabi-Yau manifolds, the complete
intersections
in toric manifolds or homogeneous spaces, this mirror symmetry prediction
\cite{Ca, BV, BCKV1, BCKV2}
can be interpreted as a quantum cohomology counterpart of the
weak Lefschetz hyperplane section theorem relating
cohomology algebras of the ambient manifolds and their hyperplane sections.
As it is mentioned in \cite{GS}, \lq\lq quantum hyperplane section conjecture"
can be formulated in intrinsic terms of Gromov-Witten          
theory on the ambient manifold and does not require a reference
to its mirror partner.
In this paper we formulate and prove the conjecture for homogeneous spaces.
It would be one of the highly nontrivial functorial properties
enjoyed by quantum cohomology algebras.  One
can compute the virtual numbers of rational curves on a Calabi-Yau
3-fold complete intersection, provided one knows the
quantum cohomology algebra of the ambient space. In fact, one needs to know the
quantum differential equations of the space, which are certain
linear differential equations arising from the flat
connection in the quantum cohomology algebra.
The mirror symmetry prediction is that
the quantum differential equations of a Calabi-Yau manifold are
equivalent (in a sense) to the Picard-Fuchs differential
equations of another Calabi-Yau manifold. In contrast,
the proposed conjecture is that
there is a certain relation between quantum differential
equations of a manifold and those of a certain complete
intersection.
When the ambient space is a symplectic toric manifold, the conjecture is
a corollary of the Givental mirror theorem \cite{GT}.
                                                                  
\bigskip

Let $X$ be a compact
homogeneous space of a semi-simple complex Lie group
and let $V$ be a vector bundle over $X$.
Suppose $V'_\beta :=\pi _*e_1^*V$ becomes a vector orbi-bundle
over Kontsevich moduli space $\overline{M}_{0,0}(X,\beta)$
where $e_1$ is the evaluation map at the (first) marked point
from $\overline{M}_{0,1}(X,\beta )$ to $X$  and
$\pi$ is the map from $\overline{M}_{0,1}(X,\beta)$
to $\overline{M}_{0,0}(X,\beta)$ associated with \lq\lq forgetting
the marked point" \cite{Ko}.
Then one might want compute
\[ \int _{\overline{M}_{0,0}(X,\beta)}Euler(V'_{\beta}). \]
Introduce a formal parameter $\hb$. Then
it turns out that the classes
\[
G_{\beta}^V
:=(e_1)_*\frac{Euler(V_{\beta})}{\hb (\hb -c)}\]
would be better considered \cite{GE}, where
$V_\beta =\pi ^*(V'_\beta )$ and $c$ (depending on $\beta$)
are the first Chern classes
of the universal cotangent line bundles.
The classes are in $H^*(X) [\hb ^{-1}]$.
They recover the original integrals which we want:
\[ \int _X G^V_\beta  =
\frac{-2}{\hb ^3}\int _{\overline{M}_{0,0}(X,\beta)}Euler(V'_{\beta})
+o(\hb ^{-3}).\]

Consider the classes  \[ G_{\beta}^X:=(e_1)_*
\frac{1}{\hb(\hb -c)}\]
corresponding to $X$ itself (without $V$).
When $V$ is a convex, decomposable, vector bundle $\oplus L_j$ of
line bundles $L_j$,
the main result of this paper proves some explicit
relationship between $A:=\{ G_{\beta}^V |\  \beta\in H_2(X,\ZZ ) \}  $
and $B:=\{ H^V_{\beta}\cup G_{\beta}^X | \ \beta\in H_2(X,\ZZ ) \}$, where
\[ H^V_\beta = \prod _j\prod _{m=0}^{<c_1(L_j) ,\beta >}
(c_1(L_j) + m\hb ) \]
which is the key object introduced in this sequel.

\bigskip
We now formulate the precise result of this paper.
Let
$\{ p_i \}_{i=1}^k$ denote the $\ZZ _+$ basis of the
closed integral K\"ahler (ample) cone of $X$.
Let us introduce formal parameters $q_i$, $i=1,...,k$,
and the ring $\QQ [[ q_1,...,q_k]]$ of formal power series of $q_i$.
Denote by $q^\beta$
\[ \prod _{i=1}^k q_i^{<p_i,\beta >}.\]
For simplicity, let $G^X_0 =1$ and $G^V_0=Euler(V)$.
We want to compare generating functions $J^V$ and $I^V$ from $A$ and $B$,
respectively:
\beqn S^V &:=& \sum _{\beta} q^\beta G^V_\beta \\
 \Phi ^V &:=& \sum _{\beta }q^\beta H_\beta ^V \cup G^X_\beta .\eeqn
We prove that one can be transformed to another by a
unique \lq\lq mirror" transformation. To describe the transformation, let
\[ q_i = e^ {t_i}, \text{ for } i=1,...,k, \]
and introduce another formal variable $t_0$. Define degree
of $q_i$ by
\[ c_1(TX )
-c_1(V)=\sum (\deg q_i) p_i .\]                          
Let
\beqn J^V (t_0,...,t_k) &:=& e^{(t_0+\sum _i p_it_i)/\hb }S^V  \eeqn
and
\beqn
 I^V (t_0,...,t_k)&:=&  e^{(t_0+\sum _i p_it_i)/\hb }\Phi ^V ,\eeqn
which are formal power series  of $t_1,...,t_k,
e^{t_0},...,e^{t_k}$ over $H^*(X)[\hb ^{-1}]$.

{\theorem\label{thmmain} Assume that $\deg q_i\ge 0$ for all $i$. Then
$J^V$ and $I^V$ coincide
up to a unique weighted homogeneous change of variables:
$t_0\mapsto t_0+f_0\hb +f_{-1}$
and $t_i\mapsto t_i +f_i$,
where $f_{-1},...,f_k$ are power
series of $q_1,...,q_k$ over $\QQ$ without constant terms,
$\deg f_i=0$, $i=0,...,k,$ and $\deg f_{-1}=1$.}

\bigskip

{\it Remarks:}
0. $J^V$ will be shown to be the cohomological expression
of solutions to quantum differential equations associated
to $(X,V)$, which is closely related to the quantum differential
equations of the smooth zero locus of $V$.
The theorem can be extended to the case of decomposable
concavex vector bundles $V$.
                                                          
1. The change of variables is uniquely determined by
coefficients of $1=(\frac{1}{\hb})^0$ and $\frac{1}{\hb}$
in the expansions of $J^V$ and $I^V$ as power
series of $\frac 1\hb$.

2. In the case of a symplectic toric manifold $X$
the similar statement is a corollary of a mirror theorem in \cite{GT},
where $\Phi ^X$ is explicitly known.

3. For the proof of \ref{thmmain} we follow the scheme of Givental's
proof \cite{GE, GT} of the mirror theorem for nonnegative complete
intersections in toric manifolds.

4. The theorem verifies the prediction \cite{BCKV1} of virtual numbers in
Calabi-Yau 3-fold complete intersections in Grassmannians.

5. A mirror construction is established for complete intersections
in partial flag manifolds \cite{BCKV1, BCKV2}.
Because of the known quantum cohomology structure \cite{Ionut},
in principle there is no essential difficulty in finding $G^X_\beta$ for
each partial flag manifold $X$, even though a general formula of it is
unknown.

6. The quantum hyperplane section principle is applied to a nonconvex
manifold in \cite{To}.

\smallskip

{\it Notation:}
$X$ will always be a generalized flag manifold $G/P$,
where $G$ is a complex semi-simple Lie group and $P$ is
a parabolic subgroup.
Let $T$ be a maximal torus of $G$ in $P$ and let $T$ act on $X$ on
the left.
Let a complex torus $T'$ act on $X$ trivially and let
$V$ be a $T\ti T'$-equivariant convex vector bundle over $X$.
Consider $E$ a multiplicative class and suppose that
$E(V)\in H_{T\ti T'}(X)$ is invertible in
$H_{(T\ti T')}(X):= H_{T\ti T'}(X)\ot H_{T\ti T'}$,
where $H_{(T\ti T')}$ is the quotient field of  $H_{T\ti T'}(pt)$.
In section 2, we will not consider $T$-action on $X$.
In section 6, additionally we will assume that $V$ is decomposable.
The convexity of $V$ is by definition that
$H^1(\PP ^1, f^*V)=0$ for any morphism $f:\PP ^1\ra X$.
Let $T\ti T'$ equivariant line bundles $U_i$, $i=1,...,k$,
form ample basis of ordinary Picard group.
We denote $\int _X ABE(V)$
by $<A,B>^V_{0}$, for $A, B\in  H^*_{(T\ti T')}(X)$
and also we use $\int _V A:=\int AE(V)$ (equivariant push forwards).
The Mori cone of $X$ will be denote by $\Lambda$, which               
can be identified with $\ZZ _+^k$ with respect to
coordinates $p_i:=c_1(U_i)$.
On the additive group $\ZZ ^{k}$
we will give the standard partial ordering, so
that $d:=(d_1,...,d_k)\ge 0$ means $d_i\ge 0$.
Let $\phi _v$ denote the equivariant pushforward
of $1$ under the embedding $i_v$ of
the fixed point $v$ to $(X,V,E)$; this $(X,V,E)$ has the
Frobenius structure by pairing $<,>^V_0$, so that
$A_v:=<A,\phi _v>_0^V=i_v^*(A)$ for $A\in H^*_{(T\ti T')}(X)$.
For a $G$-manifold $M$, let
$M^G$ denote the set of $G$-fixed points of $M$.
We will say simply degree and dimension for
complex degree and complex dimension, respectively.
Let $\sum _a T_a\ot T^a$ be the equivariant diagonal class of $(X,V,E)$
in $X\ti X$. That is, $<T_a, T^b>^V_0=\delta _{a,b}$.
In the paper we will consider various rings
$H^*_{T\ti T'}[[\hb ^{-1}]][[q]]$ formal power series ring
of $\hb ^{-1}, q$ over $H^*_{T\ti T'}$,  $H^*_{(T\ti T')}[[\hb ^{-1}]][[q]]$   
formal power series ring
of $\hb ^{-1}, q$ over $H^*_{(T\ti T')}$, and
$H^*_{(T\ti T')}(\hb )[[q]]$ formal power series ring
of $q$ over quotient field of $H^*_{T\ti T'}[\hb ]$.

\smallskip

{\it  Structure of the paper:}
In section 2, we recall a general theory of Gromov-Witten invariants
and quantum cohomology. We introduce the Givental Correlators $S^V$.
In section 3, we show that the equivariant correlators satisfy  certain
\lq\lq almost recursion relations."
In section 4, we introduce the double construction and
show that the correlators satisfy certain polynomiality
in the double construction.
In section 5, we introduce certain class $\mathcal{P}(X,V,E)$
of series of $q=(q_1,...,q_k), \hb ^{-1}$ over $H^*_{T\ti T'}(X)$,
where a \lq mirror' group acts freely and transitively.
In section 6, we introduce a modified correlator of $S^X$. 
It will also
belong to the class $\mathcal P (X,V,E)$.
The modification is given by the hypergeometric correcting Euler classes
$H^V_\beta$ according to the decomposition type of $V$.
In sections 7 and 8, we analyze the torus $T$ action on a generalized
flag manifold and its one dimensional orbits, the representations of
the section spaces of equivariant line bundles restricted to the orbits.
The analysis would be useful to find the explicit expression of $\Phi ^X$.

{\it Acknowledgments:}
I am grateful to A. Givental and Y.-P. Lee for helping me
to understand the paper \cite{GE}; and
V. Batyrev, I. Ciocan-Fontaine, B. Fulton, B. Kreu\ss ler,
E. Tj\o tta, K. Wirthm\"uller for useful discussions
on the papers \cite{GE, GT}.
Also, I would like to thank Institut Mittag-Leffler
for the financial support during
the year-long program, \lq\lq enumerative
geometry and its interactions
with theoretical physics" in 1996/1997.
My special thank goes to D. van Straten for numerous comments and
help to improve the clarity of the paper.
                                                
\section{Mirror Symmetry}\label{setup}
\subsection{The moduli space of
stable maps}
To fix notation we recall the definition of
stable maps and some elementary properties of
the moduli spaces of stable maps to $X$ \cite{Ko, FuP, BM}.
The notion of stable maps is due to M. Kontsevich.
We recommend the (survey)
paper of W. Fulton and R. Pandharipande \cite{FuP}.
   
A prestable rational curve $C$ is
a connected arithmetic genus $0$
projective curve with possibly nodes.
The curve is not necessary irreducible.
A prestable map $(f,C; x_1,...,x_n )$ is a morphism $f$ from $C$ to
$X$ with fixed ordered $n$-many marked distinct smooth points $x_i\in C$.
We will identify $(f,C; \{ x_i\} )$ with $(f',C',\{ x_i'\})$
if there is an isomorphism $h$ from $C$ to $C'$ preserving
the configuration of marked points such that $f=f'\circ h$.
A stable map $(f,C; \{ x_i\} )$ is a prestable map
with only finitely many automorphisms. 
 
Let $\overline{M}_{0,n}(X,\beta )$ be the (coarse moduli)
space of all stable maps $(f, C; \{ x_i \} _{i=1}^n)$
with the fixed homology type
$\beta = f_*([C])\in H_2(X,\ZZ )$.
Whenever it is nonempty, the moduli space
is a connected\footnote{For a proof of the connectedness see \cite{Th}.} 
 compact complex {\it orbifold} 
with complex dimension
$\dim X + <c_1(TX),\beta > + n -3$. 

\bigskip

More precisely,
locally near a stable map the moduli space has data of a quotient
of a holomorphic domain by the (finite)
group action of all automorphisms of the stable map.
In the paper \cite{FuP}  are
constructed smooth open complex domains $V$ with
finite groups $\Gamma$ which act on $V$ such that $V/\Gamma$ are
naturally glued together in the moduli space of stable maps.
Let $X\subset _i \PP ^N$, $\beta\ne 0$, and $(X,i_*(\beta ))\ne
(\PP ^1, \kline )$. Here $\kline$ denotes the line class of $H_2(\PP ^N)$.
Given a stable map $(f,C)$ (without marked points for
simplicity), choose hyperplanes $H_j$ in $\PP ^N$
satisfying that $\{ H_j\}$ gives rise to
a basis of $H^0(\PP ^N ,{\mathcal O}(1))$,
$f$ is transversal to the hyperplanes,
and their inverse images $\{ x_{i,j} \} _i=f^{-1}(H_j)$
contain no nodes of $C$.
Then the data $(C;\{ x_{i,j}\})$ determines a point in
the moduli space of marked stable curves. Conversely, a point 
in a suitable closed subvariety of an
open smooth domain of the moduli space of marked stable curves
naturally determines a stable map $f$ with the extra 
choices of elements in $(\CC ^\ti )^N$.
If $G$ is the product of the symmetric group of the
elements of the each group $\{ x_{i,j}\}_i$, then
this
$G$ has an action sending the data $(f,C; \{ x_{i,j}\})$
to another by permutations of the new marked points.
A $(\CC ^\ti )^N$ - bundle of the smooth closed subvariety is
an algebraic local chart of the moduli space of stable maps at $f$ with the
induced $G$ action.

\smallskip

{\it Example: }
Let $X=\PP ^2$ and $f$ be a stable map without marked points
such that $f$ is transversal to the hyperplanes $x=0$, $y=0$, and $z=0$.
Assume no singular points of $C$ are mapped into the
hyperplanes, and $f_*[C]=2[\text{line}]$.
Consider their inverse images (Cartier divisors),
$a_1, a_2$, $b_1, b_2$, $c_1, c_2$ in $C$. 
This information $(C; a_1,...,c_2)$ as a stable curve 
will determine $f$ uniquely with
$(\CC ^\ti )^2$ ambiguity. This $(\CC ^\ti )^2$-bundle
over some open subset of the smooth space
$\overline{M}_{0,6}$ is the local smooth 
chart. Notice that for instance,
$(C; a_2, a_1, b_1, b_2, c_2, c_1)$ gives rise to the
same $f$ up to isomorphism. 
Thus we have to take account of the quotient by the finite
group permuting the elements of sets
$\{ a_1, a_2 \}$, $\{ b_1, b_2 \}$ and
$\{ c_1, c_2\}$.

{\it Claim:} The stabilizer subgroup $G _{(C;\{ x_{i,j}\} )}$ of
$G$ is
exactly the automorphism group $Aut(f,C)$ of $(f,C)$.

{\it Proof:} We shall construct a correspondence between 
$G _{(C;\{ x_{i,j}\})}$ 
and $Aut(f,C)$. Let $g\in G _{(C;\{ x_{i,j}\} )}$
which is given by one of the suitable permutations of $x_{i,j}$. 
So, $g(C;\{ x_{i,j}\} ) = (C; \{ g(x_{i,j})\} )$.
Since the permutation does not change the
stable curve $(C;\{ x_{i,j}\} )$, 
there is an isomorphism $h$ from 
$(C;\{ x_{i,j}\})$ to $(C;\{  g(x_{i,j})\} )$. 
The isomorphism $h$ is unique since  
there is no nontrivial automorphism in the stable curve of genus
$0$.
Of course this $h$ gives rise to an automorphism of $(f,C)$.

Conversely, if $h$ is an automorphism of $(f,C)$, then
it induces an isomorphism from $(C;\{ x_{i,j}\} )$ 
to $(C;\{ g(x_{i,j})\})$ for
a unique permutation $g$ which we allow. 
Thus we established 1-1 correspondence,
which can be easily seen to be a group homomorphism.

{\it Remark:}
The action of $Aut(f,C)$ may not be effective in
general. For instance,
see $\overline{M}_{0,0}(\PP ^1,2\kline )$.

\subsection{Gromov-Witten Invariants and $QH_{(T')}^*(V)$}

There are natural morphisms on the moduli spaces,
namely, evaluation maps $e_i$
at the $i$-th marked points and forgetting-marked-point maps $\pi$:
\[\begin{CD}
 \overline{M}_{0,n+1}(X,\beta ) @> e_{n+1} >> X \\
 @V\pi VV  \\
 \overline{M}_{0,n}(X,\beta ). \end{CD} \]
If $s_i$ are the universal sections for the marked points, then
$e_i=e_{n+1}\circ s_i$ (here we assume that $\pi$ is the
forgetful map of the last marked point).
In the orbifold charts,
$\pi$ gives the universal family of stable maps as a fine 
moduli space. 

Consider, for a second homology class $\beta \ne 0$ and
an integer $n\ge 0$, the vector orbi-bundle
$V_{\beta}=\pi _*(e_{n+1}^*(V ))$.
Here $\pi$ is a flat morphism in the level of orbifold charts.
Thus indeed, $V_{\beta}$ is vector orbi-bundle
with the fiber $H^0(C, f^*(V  ))$ at
$(f, C; \{ x_i\})$.
Notice that $V_{\beta} =\pi ^*(V_{\beta})$
(it has nothing to do with marked points).

\bigskip

{\bf Notation:} for $A_i\in H_{(T')}^*(X)$,
\begin{eqnarray*}
   V_0&:=& V \\
 \overline{M}_{0,i}(X,0)&:=& X \ \ \text{ for } i=0,1,2 \\
<A_1,...,A_N>^V_{\beta} &:=& \int _{\overline{M}_{0,N}(X,\beta )}e_1^*(A_1)
\cup ....\cup e_N^*(A_N) \cup E(V_\beta ) 
\end{eqnarray*}
Then one can show that for all $\beta$
$$\sum _{\beta _1+\beta _2=\beta}\sum _a <A_1,A_2,T_a>^V<T^a,A_3,A_4>^V $$
are totally symmetric in $A_i$. This property will be 
equivalent to the associativity of the quantum cohomology of $QH_{(T')}^*(V)$
which we define in the below.

\bigskip

Let us choose a basis $\{ p_i \}_{i=1}^k$
of $H^2(X)$ by classes in the closed K\"ahler cone.

\bigskip

{\bf Notation:}
\beqn q^\beta &:=& \prod _i  q_i^{<p_i,\beta>} \\
<A_1,...,A_N>^V &:=& \sum q^\beta <A_1,....,A_N>^V_\beta \eeqn

\bigskip

The quantum multiplication $\circ$
is defined by the following simple requirement;
for $A, B ,C\in H^*_{(T')}(X)$ 
$$<A\circ B,C >_0^V = <A,B,C>^V$$
which is a formal power series of parameters $q_i$.
Thus our quantum cohomology $QH^*_{(T')}(V)$
is defined as $H^*_{(T')}(X)\ot _{\QQ }\QQ [[q_1,...,q_k]]$ with a product
structure.

\subsection{Givental's Correlators}
We review the topic after \cite{Du, GE}. 
\subsubsection{The flat connections and the
fundamental solutions}\label{fund}
Now let $q_i=e^{t_i}$ with the formal parameters $t_i$.
We have a one-parameter family
of the formal $\mathcal D$-module structures on $QH^*_{(T')}(V)$ by
giving a flat connection 
$\nabla _i= \hb\frac{\partial}{\partial t_i}-p_i\circ$ 
for any nonzero $\hb$, $i=1,...,k$. 
For the fundamental solutions we introduce $c_i\in H^*_{T'}(
\overline{M}_{0,N}(X,\beta ))$, so-called
Gravitational descendents. These $c_i$ are the first Chern classes of
the universal cotangent line bundles at the $i$-th marked 
points.  
The line bundles are, by definition,
the dual of the normal bundle of $s_i(\overline{M}_{0,N}(X,\beta ))$
in $\overline{M}_{0,N+1}(X,\beta )$.

\bigskip

{\bf Notation:}
Let $f_i(x)\in H^*_{(T')} [y][[x]]$ for indeterminant $x,y$.
Through out the paper, 
\beqn <A_1f_1(c),...,A_Nf_N(c);B>_{\beta }^X 
&:= &\int _{\overline{M}_{0,N}(X,\beta )}
 e_1^*(A_1)f_1(c_1) ... e_N^*(A_N)f_N(c_N) B \\
 <A_1f_1(c),...,A_Nf_N(c);B>_{\beta }^V 
&:=&  <A_1f_1(c),...,A_Nf_N(c);BE(V_\beta )>_{\beta}^X  \\
 <A_1f_1(c),...,A_Nf_N(c);B>^X &=& \sum _{\beta} q^{\beta}
 <A_1f_1(c),...,A_Nf_N(c);B>_{\beta }^X  \\
 <A_1f_1(c),...,A_Nf_N(c);B>^V 
&:=&  \sum _{\beta} q^{\beta}<A_1f_1(c),...,A_Nf_N(c);B>_{\beta }^V 
\eeqn
where $A_i\in H^*_{(T')}(X),$ $B\in H^*_{(T')}(\overline{M}_{0,N}(X,\beta ))$.

\bigskip

The system of the first order equations
$\nabla _i \vec{s}=0$, $i=1,...,k$, has
the following complete set of ($\dim H^*(X)$)-many
solutions \cite{GE},
\begin{eqnarray*} \vec{s}_a &:=& 
\sum _b<\frac{e^{pt/\hb}T_a}{\hb -c}, T_b>^VT^b ,
\end{eqnarray*}
where $pt$ denotes $\sum _{i=1}^k p_it_i$ and
$\hb$ is a formal variable (but when $\beta =0$, set $\hb =1$).

The following two formulas show that $\vec{s}_a$ are indeed
solutions to the quantum differential system $\nabla _i\vec{s}=0$.

Using that $c_i-\pi ^*(c_i )$ is the fundamental
class $\Delta _i$ represented by the section $s_i:
\overline{M}_{0,n}(X,\beta )\ra
\overline{M}_{0,n+1}(X,\beta )$ and $c_i\cup\Delta _i =0$ (the 
image of $s_i$ is isomorphic to 
$\overline{M}_{0,3}(X,0)\ti _X\overline{M}_{0,n}(X,\beta )$), it is easy to
derive so-called
the fundamental class axiom and the divisor axiom \cite{W, GE}.
Let $f_i(x)$ be polynomial with coefficients in 
$\pi ^*(H^*_{(T')}(\overline{M}_{0,n}(X,d)))$. Let $D$ be a divisor class
in $H^*_{(T')}(X)$. Then (for $n>0$)
\beqn
<f_1(c),...,f_n(c),1>^V_{\beta}
&=&\sum _i<f_1(c),...,\frac{f_i(c)-f_i(0)}{c},...,f_n(c)>^V_{\beta} ,
\eeqn
(where we abuse the notation \lq$f_i(\pi ^*(c))=f_i(c)$',)
and
\beqn
<f_1(c), ..., f_n(c), D>^V_{\beta} 
&=&
<D,\beta >
<f_1(c), ...,f_n(c)>^V_{\beta} \\
+
\sum _i
<f_1(c),...,f_{i-1}(c), &
D\frac{f_i(c)-f_i(0)}{c} &,f_{i+1}(c),..., f_n(c)>^V_{\beta}.
\eeqn

\bigskip

Consider
\[e^{pt/\hb }S^V :=
\sum _a<\vec{s}_{a},1>_0^VT^a=
e^{pt/\hb }\sum _a <\frac{T_a}{\hb(\hb -c)}>^VT^a=e^{pt/\hb }(1+o(1/\hb )). \]
which is the main object in this paper. This $S^V$
will be called the Givental's correlator for $(X,V,E)$.
It is an element in $H^*_{(T')}(X)[\hb ^{-1}][[q_1,...,q_k]]$.
Notice that $S^V$ for $(X,V,Euler)$ is homogeneous of degree $0$
if we let $\sum (\deg q _i)p_i= c_1(TX)-c_1(V)$,  $\deg \hb =1$,
and $\deg A=b$ if $A\in H^{2b}_{T'}(X)$.

\bigskip

The quantum $\mathcal D$-module of $QH^*_{(T')}(V)$ is defined
by the $\mathcal D$-module generated by $<\vec{s},1>^V_{0}$ for all
flat sections $\vec{s}$.
When there is no $V$ considered, we denote by $QH^*(X)$ the
quantum cohomology. That is, using $<...>^X$, we define
$QH^*(X)$.

\bigskip

{\it Remark:}
Suppose a differential operator $P(\hbp ,e^{t_i} ,\hb )$ with
coefficients in $H^*_{(T')}$
annihilates $<\vec{s},1>^V_0$ for all flat sections $\vec{s}$, 
then $P(p_1,...,p_k,q_1,...q_k,0)$
holds in $QH^*_{(T')}(V)$ \cite{GE}. 

\subsubsection{Examples}
The projective space $\PP ^{n}$:  It is
well known that  in the quantum cohomology ring $QH^*(\PP ^n)$,
$(p\circ )^{n+1}=q$, where $p=c_1({\mathcal O}(1))$
and $q$ is given with respect to the line class dual to $p$.
The corresponding operator is $(\hb\frac d{dt})^{n+1} - e^t$ .
The solutions are 
explicitly known in \cite{GH}. $S^X$ is
      \[   1+ \sum _{d>0}e^{dt}
\frac 1{((p+\hb )(p+2\hb )...(p+d\hb ))^{n+1}}.\]

The complete flag manifolds $F(n):$ Let $F(n)$ be the set of all
complete flags $(\CC ^1\subset ...\subset \CC ^n)$ in $\CC ^n$.
The usual cohomology ring is
$\QQ [x_1, x_2,..., x_n]/(I_1, ..., I_n)$
where $x_i$ are the Chern classes of $(S_i/S_{i-1})^*$,
$S_i$ are the universal subbundles with fibers $\CC ^i$  and
$I_i$ are the $i$-th elementary symmetric polynomials of
$x_1,...,x_n$.
Let us use as a basis of $H_2(F(n),\ZZ )$ duals of 
the first Chern classes of $(S_i)^*$, $i=1,...,n-1$.
They are in the edges of the closed K\"ahler cone.

Let $A(x_i)$ be a matrix
\[\left( \begin{array}{cccccc} 
          x_1 & q_1 & 0   & 0   & ...     & 0  \\
          -1  & x_2 & q_2 & 0   & ...     &    \\
              & ... &     &     & ...     &    \\
           0  & ... &     &  -1 & x_{n-1} & q_{n-1} \\
           0  & ... &     &   0 & -1      & x_n
\end{array}\right).\]
Then the quantum relations are generated 
by the coefficients of the characteristic
polynomial of the matrix $A(x_i)$.

The corresponding differential operators turn out to be
obtained by the same method using $A(x_i)$ with arguments
$\hb\frac{\partial}{
\partial t_1}$ instead of $x_1$, 
$\hb\frac{\partial}{
\partial t_{i}}-\hb\frac{\partial}{
\partial t_{i-1}}$ instead of $x_i$, and
$-\hb\frac{\partial}{
\partial t_{n-1}}$ instead of $x_n$ \cite{GS, KT}. 
These differential operators are the integrals
of the quantized Toda lattices.
The quadratic differential operator of them
can be easily derived. 
In fact, given a quantum relation of $F(n)$ between
the divisors $x_i$, there is a unique operator satisfying that 
its symbol becomes the relation and it
annihilates $<\vec{s},1>_0$ for all flat sections $\vec{s}$.
In general, the explicit cohomological expression $S^X$ of
solutions to the quantum differential operators are not known.

\subsubsection{The general quintic hypersurface
in $\PP ^4$} Let $Y$ be
a smooth degree $5$ hypersurface in $\PP ^4$.
$Y$ is not a homogeneous space. However, using virtual
fundamental class $[\overline{M}_{0,n}(Y,\beta )]$ \cite{Be, BF,
LT}, one can
define also the quantum cohomology $QH^*(Y)$ of $Y$.
It is expected that
\[ <A_1f_1(c),...,A_Nf_N(c)>^Y = <A_1f_1(c),...,A_Nf_N(c)>^{{\mathcal O}(5)}. \]
Let $p$ be the induced class of
the hyperplane divisor in $\PP ^4$. The quantum
relation is $(p\circ ) ^4=0$.
The corresponding operator is, 
however, {\it not}  $(\hb\frac d{dt})^4$,
but  $$(\hb\frac d{dt})^2
\left( \frac{(\hb\frac d{dt})^2}{5+f(q)}\right),$$ where
$<p\circ p,p>^Y_{0} = 5+f(q)$.
Notice that in this $Calabi-Yau$ 3-fold case,
we lose the whole information of quantum cohomology when one
concerns only the quantum relation, $(p\circ )^4=0$. 
The unknown $f(q)$ was conjectured by physicists
\cite{Ca}. The general idea of the prediction 
is the following. Roughly speaking,
in theoretical physics, there are quantum
field theories
associated to Calabi-Yau three-folds by $A$-model and $B$-model.
What we have constructed so far are $A$-model objects for
Calabi-Yau three-folds. On the other hand using a
family of the so-called mirror manifolds which are also Calabi-Yau
three-folds, conjecturally
one may construct the equivalent quantum field theory by $B$-model.
The corresponding mirror partner of a  quantum differential 
equation / quantum ${\mathcal D}$ module 
is the Picard-Fuchs differential equation 
/ Gauss-Manin connection of the
mirror family. It was predicted that
they are equivalent by a certain transformation.
In \cite{Ca} are obtained the conjectural mirror family of quintics,
the Picard-Fuchs differential equation and the transformation.
That is how the prediction is made. The prediction is now proven to be correct
by Givental \cite{GE}.

\subsection{The idea of the proof of theorem \ref{thmmain}}
To describe the idea, let us
notice that Givental's proof \cite{GT} of the mirror
conjecture for the nonnegative toric complete intersections
can be divided into three parts. (He shows in the paper that
the mirror phenomenon occurs also in non-Calabi-Yau manifolds.)
 Let $X$ be a Fano toric manifold
with a big torus $T$,
and $V$ be a $T\ti T'$-equivariant decomposable convex vector bundle over $X$,
where $T'$ acts on $X$ trivially.

\begin{enumerate}
\item
In $A$-part, it is proven that
\begin{enumerate}
\item
 the $T\ti T'$-equivariant
solution vector $S^V\in H^*_{(T\ti T')}(X)[[q,\hb ^{-1}]]$ has an
 \lq\lq almost recursion relation,"
\item
 it satisfies the polynomiality in the so-called
\lq\lq double construction," and
\item
 it is uniquely determined by the above two properties with
the aymptotical behavior $S^V=1+o(\frac 1{\hb})$.
 
\end{enumerate}
 
\item In $B$-part, another
($T\ti T'$-equivariant hypergeometric) vector $\Phi ^V$, presumably
given by the mirror symmetry conjecture, is constructed.
It is verified that it also satisfies (a) and (b)
using a toric (naive) compactification of holomorphic maps from
$\PP ^1$ to $X$.
 
\item
When $c_1(X)-c_1(V)$ is nonnegative and $E=Euler$,
there is a suitable equivalence transformation between
$\Phi ^V$ and $S^V$.
 
\end{enumerate}
 
In this paper, for a $T\ti T'$ equivariant decomposable
convex vector bundle $V$ over
any compact homogeneous $X$ of a semi-simple complex
Lie group $G$, we will show that $S^V$ 
satisfies property 1 above. In this case, $T$ is 
a maximal torus of $G$.
 
We define $\Phi ^V$ which corresponds $\Phi ^V$ of the toric case
in property 2:
Let $\Phi ^X=S^X=\sum _d\Phi ^X_{d}q^d$.
For $\Phi ^V$, we will find
a modification $H'_d\in H^*_{T\ti T'}(X)[\hb ]$ 
(depend on $V$ and $d$) such that
if $\Phi ^V :=\sum _d \Phi ^X_{d}H'_dq^d$,
then
 
(A) $\Phi ^V$ (after the restriction to
 the fixed points) has the almost recursion relation exactly like $S^V$
and
 
(B) $\Phi ^V$ has 
the polynomial property in the double construction.
 
In fact, we design $H'_d$ to satisfy (A) and (B).

Finally, when $E=Euler$ and $c_1(TX)-c_1(V)$ is nonnegative,
we will prove that a certain operation
will transform $S^V$ to $\Phi ^V$, since they satisfy
the same almost
recursion relation and the polynomiality of the double construction.

\section{The almost recursion relations}\label{sectionre}
As in section \ref{setup}
let $X$ be a homogeneous manifold $G/P$ where $G$ is a complex semi-simple
Lie group and $P$ is a parabolic subgroup. Let $T$ be a maximal torus.
The $T$ action has only isolated fixed points $\{ v,w,... \}$.
The one dimensional invariant orbit of $T$ is analyzed in detail in
section \ref{grass} and \ref{flag}. For a moment we need the fact 
that the closures of orbits are {\em finite}
$\PP ^1$'s connecting a fixed point $v$ to another
fixed point $w$.
For a given equivariant vector bundle $W$ 
over a $T\ti T'$-space $M$, we use $[W]$ which denote the element in
the $K$-group
$K^0_{T\ti T'}(M)$ corresponding to the $T\ti T'$ vector bundle $W$.
\bigskip

The torus action on $X$ induces the natural action
on the moduli space of stable maps by the functorial property. 
Since the evaluation maps are $T\ti T'$-equivariant, 
the pullbacks of $T\ti T'$-bundles have natural actions 
in the orbifold sense. 
In turn, $V_\beta$ has
the induced $T\ti T'$-action. 
{\it All ingredients in section 2
are from now on the equivariant ones.}
We would like to evaluate $S^V$ as a specialization
of the equivariant one corresponding to $S^V$.
We use the same notation 
$S^V\in H^*_{(T\ti T')}(X)[[\hb ^{-1}]][[q_1,...,q_k]]$ 
for the equivariant one.
Notice that $S^V$ might have power series of $\hb ^{-1}$ in each coefficient
of $q^d$, since $c$ are not anymore nilpotent.
Using the localization theorem,
we shall find an \lq\lq almost recursion relation" on the equivariant
Givental correlator.
To begin with,
we summarize the fixed points of the induced action on
the moduli space of stable maps.
If a stable map represents a fixed point in the moduli
space, the image of the map should lie in the
closure of the 1-dimensional
orbits. The special points are mapped to isolated fixed points.
Let us denote by $\ka _{v,w}$
the character of the tangent space of a
1-dimensional orbit connecting an isolated fixed point $v$ to another $w$.
Then, $-\ka _{v, w}$ is the character of the tangent
line of the 1-dimensional
orbit $o(v,w)$ at $w$. We use $\beta _{v,w}$  to stand for
the second homology class represented by the ray.
Denote by $o(v)$ the set of all fixed points $w\ne v$ which
can be connected by a one-dimensional orbit with $v$.      

\bigskip

{\lemma\label{lemmare} {\bf Recursion Lemma} {\em (\cite {GE})}

Denote by $\phi _{v}$ the equivariant classes $i_*(1)$ at $v$,
here $i_v$ denotes the $T\ti T'$-equivariant
inclusion of the point $v$ into $(X,V)$.
Then $S^V\in H^*_{T\ti T'}(X)[[\hb ^{-1}]][[q_1,...,q_k]]$ 
has an \lq\lq almost" recursion relation, namely, for any $v\in X^{T}$,

0) $S^V_v(q,\hb ):=<S^V,\phi _v>^V_0 \in H^*_{(T\ti T')}(\hb )[[q]]$ 
and the substitution
$S_w(q,-\ka _{v,w}/m)$ of $\hb$ with $-\ka _{v,w}$ in $S_w(q,\hb )$
 is well-defined,
 
1) The difference $R_v$ of $S^V_v(q,\hb )$ and the \lq\lq recursion part"
is a power series of $q$
over the {\em polynomial} ring of $1/\hb$, that is,
 \[ R_v:=
S^V_v(q,\hb ) -
\sum _{w\in o(v),\ m>0} q^{m\beta _{v,w}}
\frac {(-\ka _{v,w})/m}{\hb (\ka _{v,w}+m\hb)} 
\frac{E(V_{v, w,m})i_v^*(\phi _v)}
{Euler(N_{v,w,m})}S_w(q,-\ka _{v,w}/m),
\]
is in $H^*_{(T\ti T')}[\hb ^{-1}][[q]]$, where
$V_{v,w,m}$ is $T\ti T'$ representation space $H^0(\PP ^1,f^*V)$,
here $f$ is the totally ramified $m$-fold map onto $o(v,w)$
over $v$ and $w$; and
$N_{v,w,m}$ is the $T\ti T'$-representation space $[H^0(\PP ^1,f^*TX)]-[0]$;
and

2)  furthermore, for $S^X$ itself, the first term $R_v$
is $1$.}

\bigskip

We will say that the statement 1) reveals the almost recursion relation
of $S^V$.
The statement 2) shows that $S_v^X$ have  recursion relations in the
ordinary sense.

\bigskip
{\bf Proof.} 
First of all, using the short exact sequence
\[
 0 \ra Ker \ra V_d \ra e_1^*(V) \ra 0  \]
over $\overline{M}_{0,1}(X,d)$, we see that
$S^V$ is indeed in $H^*_{T\ti T'}[[\hb ^{-1}]][[q_1,...,q_k]]$.
(The last map in the sequence is given by the evaluation
of global sections at the marked point.)

A connected component of the $T$-fixed loci of
the moduli space $X_d:=\overline{M}_{0,1}(X,d)$ is isomorphic to a product of
Deline-Mumford spaces with marked points from
the special points 
of the inverse image $f^{-1}(v)$ of the generic $f$
in the component for all $v\in X^{T}$.
Now fix a $v$ and consider $S^V_v$.
It is enough to count the fixed locus $F^{d,v}$ where
the marked point $x$ should be mapped to the
fixed point $v$ since $\phi _v$ can be supported only
near the point. For a stable map
$(f,C;x)$ denote by $C_1$ 
the irreducible component of $C$ containing the marked
point $x$.
Then $F^{d,v}$ is the disjoint union  of
\[ F^{d,v}_1:=\{ (f,C;x)\in F^{d,v} \ | \ f(C_1)=v \}\]
and
\[ F^{d,v}_2:=
\bigcup _{w\in o(v), m=1,...,m\beta _{v,w}\le d} F^{d,v,w,m}, \]
where $F^{d,v,w,m}$ is
\[ \{ (f,C;x)\in F^{d,v} \ | \ 
w\in f(C_1),\ \deg f|_{C_1}=m \}. \]
 
$S^V_v$ is an integral over $F^{d,v}$'s  by a
localization theorem for orbifolds.
We claim that the integral of 
\[ \frac{E(V_d)e_1^*(\phi _v)}{\hb (\hb -c )} \]
over $F^{d,v}_1$ is in $H^*_{(T\ti T')}[\hb ^{-1}]$.
The reason is that 
the universal cotangent line bundle over $F^{d,v}_1$
in the moduli space has the trivial action. It implies that
the equivariant class $c$ restricted to $F^{d,v}_1$ is nilpotent.
 
\smallskip

Now we shall obtain the \lq almost recursion relation' 
from the contribution of the fixed loci $F^{d,v}_2$.
Denote $d-m\beta {v,w}$ by $d'$.
Since $C_1$ is always one end of $C$ for any $(f,C;x)\in F^{d,v}_2$,
we can have a natural isomorphism from
$F^{d',w}$ to $F^{d,v,w,m}$, where $F^{d',w}$ are
fixed loci in $X_{d'}:=
\overline{M}_{0,1}(X,d')$, consisting
of the stable maps sending the marked points to $w$.
We obtain the morphism, joining the $m$-covering of $o(v,w)$ to stable
maps in $F^{d',w}$. By the $m$-covering of $o(v,w)$, we
mean a totally $m$-ramified map from $\PP ^1\cong C_1$ 
to $o(v,w)$ over $v$ and $w$. Let $x'=f^{-1}(w)\cap C_1$.

\bigskip

We claim that the normal bundles as in 
$K^0(F^{d, v,w,m}\cong F^{d',w})$ satisfy  the equality
\begin{eqnarray}
[N_{X_d/F^{d, v,w,m}}]-[N_{X_{d'}/F^{d',w}}] &=&
[N_{v,w,m}] -[T_wX]  
+[T_{x'}C_1\ot L|_{F^{d',w}}]  \label{normal}
 \end{eqnarray}
where $L$ is the universal tangent line bundle 
over $X_{d'}$.
The reason of the claim is as follows:
Recall that each fixed component is isomorphic to the product
of moduli space of stable curves (see section 3 in  \cite{Ko} for detail). 
Hence, we conclude that
$[N_{X_d/F^{d, v,w,m}}]-[N_{X_{d'}/F^{d',w}}]
-[L|_{F^{d',w}}]$ (over each
fixed components)
is equal to a trivial bundle with nontrivial actions.
The twister by action can be computed by study of action
on normal spaces at $(f,C_1\cup C_2 ; x)\in F^{d,v,w,m}$.
Let $N_1$ be the normal space of $F^{d,v,w,m}$ at $(f,C_1\cup C_2 ;x )$
and $N_2$ be the normal space of $F^{d',w}$
at $(f|_{C_2}, C_2; x':=C_1\cap C_2)$. Then
as representation spaces 
\beqn [N_1] &= &[N_2]+([H^0(C_1, f|_{C_1}^*TX)]-[H^0(C_1,TC_1)]) -[T_{w}X] \\
&+& [T_{x'}C_1\ot T_{x'}C_2] +[T_{x'}C_1]+[T_xC_1].\eeqn
Hence we conclude the claim (\ref{normal}) after canceling of $[H^0(C_1,TC_1)]
=[0]+[T_{x'}C_1]+[T_xC_1]$.

\bigskip

On the other hand,
the direct sum of
the fiber of $V_d$ at $(f,C_1\cup C_2 ; x)\in F^{d,v,w,m}$
and $V|_w$
is equal to the direct sum of the fiber of
$V_{d'}$ at $(f|_{C_2}, C_2;x')$
and $H^0(C_1, (f|_{C_1})^*V)$. 

Thus, applying the localization theorem we obtain
\beqn
&&\int _{X_d}\frac{E(V_d)e_1^*(\phi _v)}{\hb(\hb -c)}=I+ 
\sum _{w\in o(m), 0<m; m\beta _{v,w}\le d}
\frac{E(V_{v,w,m})i^*_v(\phi _v)(-\ka _{v,w}/m)}{m\hb (\ka _{v,w}/m +\hb )
Euler(N_{v,w,m})} \\
&& \ti \int _{X_{d-m\beta _{v,w}}}
\frac{E(V_{d-m\beta _{v,w}})e_1^*(\phi _w)}{(-\ka _{v,w}/m)(-\ka _{v,w}/m - c)},
\eeqn
where $I$ is the integral over $F^{d,v}_1$.
The factor $m$ in $m\hb (\ka _{v,w}/m +\hb ) $ comes from the nature of 
orbifold localization theorem. (There are $m$ automorphisms of $f|_{C_1}$.)

\bigskip

Using induction on $|d|=\sum d_i$,
we may assume that the integral factors in the second term
are  well-defined and belong to $H^*_{(T\ti T')}$.
(Localization theorem itself also explains them.)
So, statements 0) and 1) in the lemma are proven.

\smallskip

Now let us prove statement 2).
Since
$<c_1(TX),\beta >\ge 2$ for all $\beta$, by
degree counting we see that there are no contributions from
the integral over $F^{d,v}_1$. The reason is that
$\dim \overline{M}_{0,\sum d_i+1} =(\sum d_i)-2$ is less than
$2(\sum d_i )-2$ if $(d_1,..,d_k)\ne 0$ and 
$\dim \overline{M}_{0,1}(X,d)\ge 2\sum d_i+\dim X -2$.
So, in the case of $S^X$, $R_v =1$.

\section{The double construction}
{\lemma\label{lemmado}{\bf Double Construction Lemma}

 The double construction
\[ W(S^V):=\int _{V} S^V (qe^{\hb z},\hb )e^{\sum p_iz_i}S^V (q,-\hb )\]
is a power series of $q_1,...,q_k$ and $z_1,...,z_k$ with 
coefficients in $H_{T\ti T'}^*[\hb ]$.}

\bigskip

A priori $W(S^V)$ has coefficients in Laurent power series ring of
$\hb ^{-1}$ over $H_{T\ti T'}^*$.
For the proof we will make use of graph spaces 
and universal classes defined in the below.

\subsection{The main lemma}
Let $L_d$ be the projective space of the  
collection of all $(f_0,...,f_N)$ such that
$f_i(z_0,z_1)$ are homogeneous polynomials of degree $d$.
$L_d$ is isomorphic to $\PP ^{(d+1)(N+1)-1}$.
Given a stable map of degree $(d,1)$ from a prestable curve
$C$ to $\PP ^N\ti \PP ^1$, 
there is a special irreducible component $C_0$ of $C$ such that $C_0$
has degree $(d_0,1)$ under the stable map. This special
component $C_0$ is parameterized by $\PP ^1$ in the target space.
Thus we can identify $C_0$ with $\PP ^1$ and keep track where 
the other components intersect. Suppose the other
connected components $C_1,...,C_l$ of $C-C_0$ 
intersect with $C_0=\PP ^1$ at $[x_1: y_1],...,[x_l: y_l]$.
If the degrees of $C_i$ are $d_i$ under the stable map,
we now associate the stable map to 
$$\prod _{i=1}^l(y_iz_0 -x_iz_1)^{d_i}(f^0_0,...,f_N^0),$$ 
where $(f^0_0,...,f_N^0)$ are the polynomials coming from
the data of the restriction of $f$ to $C_0$.

\bigskip

{\bf Main lemma: } (Givental \cite{GE})
{\em  The above \lq\lq polynomial"
mapping from $G_d(\PP ^N):=\overline{M}_{0,0}(
\PP ^N\ti \PP ^1, (d,1))$ to $L_d$ is a 
$(\CC ^\ti)^N\ti \CC ^\ti$-equivariant morphism,
where $\PP ^N$ has the diagonal $(\CC ^\ti)^N$ action
and $\PP ^1$ has the $\CC ^\ti$ action by $[z_0:z_1]\mapsto
[tz_0:z_1]$ for $t\in \CC ^\ti$.}

\bigskip

Notice that the $\CC^\ti$ action on $L_d$ is given by
\[ [f_0(z_0,z_1):...:f_N(z_0,z_1)]
\mapsto [f_0(t^{-1}z_0,z_1):...:f_N(t^{-1}z_0,z_1)]\]
for $t\in \CC ^\ti$.

\subsection{The universal class} 
The $T\ti T'$-equivariant spanned line bundle $U_i$ over $X$
gives rise to the $T\ti T'$-equivariant morphism $\mu ^i_0: X\ra\PP ^N$,
and
so we obtain:
\[
\begin{CD}
(\mu _d^i) ^*({\mathcal O}(1)) @. @.  {\mathcal O}(1) \\
@VVV  @.                 @VVV \\ 
G_d(X)@>>> G_{d_i}(\PP ^N)@>>> L_{d_i} \\
\end{CD},
\]
where $G_d(X)$ is the graph 
space $\overline{M}_{0,0}(X\ti \PP ^1, (d,1))$,
and $\mu _d^i$ is the $T\ti T'\ti \CC ^\ti$-equivariant
map from $G_d(X)$ to $L_{d_i}$.
On ${\mathcal O}(1)$ we choose the lifted $\CC ^{\ti}$-action coming from 
the action on the vector space of $N+1$ $d_i$-homogeneous polynomials by
\[ [f_0(z_0,z_1):...:f_N(z_0,z_1)]
\mapsto [f_0(z_0,tz_1):...:f_N(z_0,tz_1)]\]
for $t\in \CC ^\ti$. 

Denote by $P_i=c_1((\mu ^i_d) ^*{\mathcal O}(1))$, the $T\ti T'\ti
\CC ^\ti $-equivariant
Chern class. It is said to be
a universal class in the paper \cite{GT}.
   
Denote by $W_d$ the vector orbi-bundle
over $G_d(X)$ with the fiber $H^0(C,\psi ^*\pi _1^*V)$
at $(C,\psi )$: Consider
\[\begin{CD}
G_{d,1}(X)@>>{e_1}> X\ti \PP ^1 \\
@V\pi VV   @AAA \\
G_d(X) @. \pi _1^*V,
\end{CD}\]
where $G_{d,1}(X)$ denotes the graphs space with one marked point
and $\pi _1$ is the projection of $X\ti \PP ^1$ to
the first factor $X$. Then $W_d:=\pi _* e_1^*\pi _1^*V$.

\subsection{Proof of Lemma \ref{lemmado}}\label{pfdo}
It is enough to show the equality
\[\sum _dq^d\int _{G_d(X)}e^{Pz}E(W_d)
=\int _{V} S^V (q,\hb )e^{pz}S^V (qe^{-\hb z},-\hb ).\]
The left integral is a $T\ti T'\ti \CC ^{\ti}$-equivariant
push forward with $\hb$ as $c_1({\mathcal O}(1))$
over $\PP ^{\infty}$ and the right one is a $T\ti T'$-equivariant push forward
with a formal variable $\hb$.   

\bigskip

We will apply localization theorem.
Let us analyze the $\CC ^\ti$-action fixed loci 
$G_d(X)^{\CC ^\ti}$ of $G_d(X)$.
$G_d(X)^{\CC ^\ti}$ is isomorphic to $\sum _{d^{(1)}+d^{(2)}=d}
\overline{M}_{0,1}(X,d^{(1)})\ti _X\overline{M}_{0,1}(X,d^{(2)})$.

\bigskip

Suppose $|d^{(1)}|+|d^{(2)}|\ne 0$. The normal bundle is as follows:

When $|d^{(1)}||d^{(2)}|=0$: The codimension is 2 (one from
the nodal condition and the other from the condition 
of the image of the nodal point).
Then the Euler class of the normal bundle is $\hb (\hb -c_0)$,
or $-\hb (-\hb -c_\infty)$, where $c_0$ and $c_\infty$ are
the Chern classes of universal
cotangent line bundles of the first marked point
over $\overline{M}_{0,1}(X, d^{(1)})$ and 
$\overline{M}_{0,1}(X, d^{(2)})$, respectively.
Here we assume the following convention:
$0=[0:1], \infty =[1:0]$, the associated equivariant line bundle
to the character 1 of the group $\CC ^\ti$ has $\hb$ as its equivariant
Chern class.
                 
When $|d^{(1)}||d^{(2)}|\ne 0$: The codimension is $4$ and
the Euler class is $\hb (\hb -c_0)(-\hb )(-\hb -c_\infty )$.
Here, for instance, $c_0\in
H^2(\overline{M}_{0,1}(X,d^{(1)})\ti _X\overline{M}_{0,1}(X,d^{(2)}))$
is the pull-back of
the Chern class of the universal cotangent line bundle
of the first factor of
$\overline{M}_{0,1}(X,d^{(1)})\ti _X\overline{M}_{0,1}(X,d^{(2)})$.

\bigskip

Let us analyze $P_i$ restricted to $G_d(X)^{\CC ^\ti}$.
Consider the commutative diagram,
\[\begin{CD}
G_{d^{(1)},d^{(2)}}(X):=\overline{M}_{0,1}(X,d^{(1)})\ti _X 
\overline{M}_{0,1}(X,d^{(2)})
@>>\mu _d^i > L_{d_i} @. 
\ni z_0^{d^{(1)}_i}z_1^{d^{2)}_i}[x_0:...:x_N]\\
@V{\pi _2} VV @AAA @AAA\\
\overline{M}_{0,1}(X,d^{(2)}) 
@>>{\mu ^i_0\circ e_1}> \PP ^N @. \ni [x_0:...:x_N],
\end{CD}\]
where the first vertical map $\pi _2$ is the projection and
under the second vertical map
$\PP ^N$ is embedded into $L_{d_i}$ as the $\CC ^\ti$
-action fixed locus of the part
$\{ z_0^{d^{(1)}}z_1^{d^{(2)}}[x_0:...:x_N] | [x_0:...:x_N]\in
\PP ^N\} $.
One concludes that 
$e_1^*\circ (\mu _0^i)^*
(c_1({\mathcal O}(1)|_{\PP ^N}))=
e_1^*(p_i)-d_i^{(2)}\hb $ and so
\[ \sum P_iz_i|_{G_{d^{(1)},d^{(2)}}} 
=\sum (\pi _2^*e_1^*(p_i) -d^{(2)}_i\hb )z_i. \]

Since
\begin{eqnarray*}
&& S^V(q,\hb )e^{pz} S^V(qe^{-\hb z},-\hb) \\
&=& \sum _{a ,b,d^{(1)},d^{(2)} } <\frac{T_{a}}{\hb (\hb -c)}>^V_{d^{(1)}}
              T^{a}q^{d^{(1)}}
    e^{pz}<\frac{T^{b}}{-\hb (-\hb -c)}>^V_{d^{(2)}}T_{b}q^{d^{(2)}}e^{-d^{(2)}\hb z},
\end{eqnarray*}
we see that
\begin{eqnarray*}
& & \int _{V} S^V(q,\hb)e^{pz}S^V(qe^{-\hb z},-\hb) \\
&=& \sum _{a ,d^{(1)},d^{(2)}}<\frac{T_{a}}{\hb (\hb -c)}>^V_{d^{(1)}}
  <\frac{T^{a}e^{pz-d^{(2)}\hb z}}{-\hb (-\hb -c)}>^V_{d^{(2)}}q^{d^{(1)}+d^{(2)}} \\
&=& \sum _d q^d\int _{G_{d^{(1)},d^{(2)}}(X)}
\frac{e^{(\pi _2^*e_1^*p-d^{(2)}\hb)z}E(W_d)}
         {[N_{G_d(X)/G_{d^{(1)},d^{(2)}}(X)}]} \\
&=& \sum _d  q^d\int _{G_d(X)}e^{Pz}E(W_d),
\end{eqnarray*}
after applying the localization theorem only for $\CC ^\ti$ action
on $G_d(X)$.

\section{The class ${\mathcal P}({\mathcal C})$ 
and mirror transformations}\label{sectiontr}
\subsection{The class $\mathcal P({\mathcal C})$}
Let
${\mathcal C}$ be the collection of
given data of 
$C_{v,w,m}\in H^*_{(T\ti T')}$, $\ka _{v,w}\in H^*_{T\ti T'}$, and
$\beta _{v,w}\in \Lambda -0$, for all
$(v,w,m)\in X^{T}\ti X^{T}\ti \NN$ with $v\in o(w)$. Here $\NN$ is 
the set of positive integers.
Assume that $(p_i)_w-(p_i)_v=-<p_i,\beta _{v,w}>\ka _{v,w}$ for all $i=1,...,k$.
Define degree of $\hb$ as 1. Let $q_1,...,q_k$ be formal
parameters with 
some given nonnegative degrees. Define the degree of a homogeneous
class of $H_{T\ti T'}^b(X)$ as $b/2$.
Let $\mathcal P({\mathcal C})$ be the class of all
$Z(q,\hb )\in H^*_{T\ti T'}(X)[[\hb ^{-1}, q]]$ of
homogeneous degree 0 such that
\begin{description}
\item[a)] $Z(0,\hb )=1$, 
$Z_v(q,\hb ):=<Z,\phi _v>^V_0$ (this is not depend on $V$) 
is in $H^*_{(T\ti T')}(\hb )[[q]]$ for any fixed
point $v$, and $Z_w(q,-\ka _{v,w}/m )$ are well-defined
for all $v\in o(w),m>0$ ($m$ are 
positive integers),
\item[b)] the almost recursion relation for each fixed point
$v$ holds, that is by definition, 
\beqn R_v:= Z_v(q, \hb ) -
 \sum _{m>0,w\in o(v)}
q^{m\beta _{v,w}}\frac{C_{v,w,m}}{\hb (\ka _{v,w} +m\hb )}Z_w(q,-\ka _{v,w}/m),
\eeqn
is in $H^*_{(T\ti T')}[\hb ^{-1}][[q]]$, where
 $$ q^{m\beta _{v,w}}:=\prod _i q_i^{m<p_i,\beta_{v,w}>}$$; and
\item[c)] in the double construction
\[ W(Z)(q,z):=\int _{V}Z(qe^{\hb z}, \hb )e^{\sum p_iz_i}Z(q,-\hb ), \]
is in $H^*_{T\ti T'}[\hb][[q,z]]$.
(We use the multi-index notation for $z=(z_1,...,z_k)$ and $q=(q_1,...,q_k)$.)

\end{description}

\bigskip

Whenever the data $\mathcal C$ comes from $(X,V,E)$ as in lemma
\ref{lemmare}, we denote the class by ${\mathcal P}(X,V,E)$.
So, in the case 
\[ C_{v,w,m} :=C_{v,w,m}^V:= \frac{(-\ka _{v,w})/m\ E(V_{v, w,m})i_v^*(\phi _v)}
{Euler(N_{v,w,m})},
\]
$\ka _{v,w}$ is the character of $T_vo(v,w)$,
$\beta _{v,w}=[o(v,w)]\in H_2(X,\ZZ )$,
and $$c_1(TX)-c_1(V)=\sum _{i=1,...,k}(\deg q_i) p_i.$$
So far, we proved that $S^V$ for $E=Euler$ is in class 
${\mathcal P}(X,V,Euler)$.

\bigskip

In the below we introduce 
on ${\mathcal P}({\mathcal C})$
a transformation group generated by the following three types of operations.

\begin{description}
\item[1) Multiplication by $f(q)$]
Let $f(q)=\sum _{d\ge 0}f_dq^d$, where $f_d\in \QQ$,
$f(q)$ is homogeneous of degree $0$, and $f(0)=1$.
Then $f(q)Z\in {\mathcal P}({\mathcal C})$.

\item[2) Multiplication by $\mathrm{exp} (f(q)/\hb )$]
Let $f(q)=\sum _{d>0}f_d  q^d$,
where $f_d$ are in $H^*_{T\ti T'}$.
Suppose $\deg (f(q))=1$.
Then $Z^{new}:=\exp (f(q)/\hb  )Z$ is still in $\mathcal P({\mathcal C})$.

\item[3) Coordinate changes]
Consider a transformation: $$Z\ra Z^{new}:=
\exp (\sum _if_i(q)p_i/\hb )Z(q\exp (f(q)), \hb )$$
where $f_i(q)=\sum _{d>0} f^{(d)}_iq^d$ of homogeneous degree $0$,
$f_i^{(d)}\in \QQ$, and $q\exp (f(q))=(q_1\exp (f_1(q)),...,q_k\exp (f_k(q)))$.
Then $Z^{new}$ is still in $\mathcal P({\mathcal C})$.

\end{description}

Let us call the transformation group the mirror group.

{\theorem\label{thmtr} {(\em \cite{GT})}
Suppose $\deg q $ are nonnegative
and there is at least one element
of form $1+o(\hb ^{-1} )$ in the class ${\mathcal P}({\mathcal C})$.
Then the mirror group action on ${\mathcal P}({\mathcal C})$ is free and
transitive.}

\bigskip

First, we will check 1), 2) and 3); and prove the so-called 
uniqueness lemma and then theorem above.

{\it Proof of 1).}
First, $Z^{new}:=fZ$ is homogeneous of degree $0$, $f(0)Z(0,\hb )=1$,
 and $fZ_v$ are in 
$H_{(T\ti T')}^*(\hb )[[q]]$, and of course $Z^{new}_w(q,-\ka _{v,w}/m)$ are 
well-defined.
Second,
\beqn Z_v^{new} &=&
 f(q)R_v  
+ \sum q^{m\beta _{v,w}}\frac{C_{v,w,m}}{\hb (\ka _{v,w} +m\hb )}
Z^{new}_w(q,-\ka _{v,w}/m).
\eeqn
Thus $fZ$ has the almost recursion relation.

Finally,
\beqn W^{new} 
&:=&\int _{V} Z^{new}(qe^{\hb z},\hb )e^{pz}Z^{new}(q,-\hb ) \\
&=& f(qe^{\hb z})f(q)W,
\eeqn
which still has the polynomial coefficients in $H^*_{T\ti T'}[\hb ]$.

{\it Proof of 2).}
The new $Z^{new}$ is homogeneous of degree $0$, $Z^{new}(0,\hb)=1$,
$Z^{new}_v$ are in $H^*_{(T\ti T')}(\hb )[[q]]$,
and $Z^{new}_w(q,-\ka _{v,w}/m)$ are well-defined. 
Since $\exp (\frac{f(q)}{\hb }+\frac{mf(q )}{\ka _{v,w}})
=1+(\ka _{v,w}+m\hb ) g_{\ka _{v,w},m}$ and
$g_{\ka _{v,w},m}$ is a 
$q$-series with polynomial coefficients in $H^*_{(T\ti T')}[\hb ^{-1}]$,
$Z^{new}$ has the almost recursion relation. 

Once again,
\beqn
W^{new}
&=& \exp(\frac {1}{\hb}  (f(qe^{\hb z})-f(q)))W .\eeqn
But $f(qe^{\hb z} )-f(q ) =\sum_{d>0} f_d((e^{\hb z})^d -1) q^d$ 
is a $(z,q)$-series with polynomial coefficients 
in $\hb H^*_{T\ti T'}[\hb ]$.

{\it Proof of 3).}
The $Z^{new}$ is homogeneous of degree $0$, $Z^{new}(0,\hb )=1$,
$Z^{new}_v$
are in $H^*_{(T\ti T')}( \hb )[[q]]$, and $Z_w^{new}(q,-\ka _{v,w}/m)$
make sense.

Since $(p_i)_w-(p_i)_v=-<p_i,\beta _{v,w}>\ka _{v,w}$,
\beqn && \sum _if_i(q)(p_i)_v/\hb    
= \sum _i f_i(q)(p_i)_w/(-\ka _{v,w}/m)  \\
&& - m\sum _i <p_i,\beta _{v,w}>f_i(q)  
+ \sum _i\frac{f_i(q)(p_i)_v }{\ka _{v,w}\hb }(m\hb  +\ka _{v,w}) .
\eeqn
The exponential of the last term on the right
can be denoted by $1+(\ka _{v,w}+m\hb )g_{\ka _{v,w},m}$
where $g_{\ka _{v,w},m}$ is a $q$-series with coefficients which
are in $H^*_{(T\ti T')}[\hb ^{-1}]$.
$Z^{new}$ satisfies the almost recursion relation.

Consider the double construction
\beqn
W^{new}(q,z) &=&
\int _{V}e^{f(qe^{\hb z})p/\hb } Z(qe^{\hb z}e^{f(qe^{\hb z})},\hb )
e^{pz}e^{-f(q)p/\hb }Z(qe^{f(q)},-\hb ) \\
&=& W(qe^{f(q)}, z+\frac{f(qe^{\hb z})-f(q)}{\hb }) . \eeqn

But since $f(qe^{\hb z})-f(q)$ is divisible by $\hb $, $W^{new}$
is a polynomial $(q,z)$-series.

\bigskip

{\lemma\label{lemmaun} {\bf Uniqueness Lemma}

Let $Z=\sum _{d\ge 0}Z_dq^d$ and $Z'=\sum _{d\ge 0} Z'_dq^d$ be series in 
${\mathcal P}({\mathcal C})$.
Suppose $Z\equiv Z'$ modulo $(\frac 1\hb )^2$.
Then $Z'=Z$.}

\bigskip

{\bf Proof.}
We may suppose that
$Z'_{d}=Z_{d}$ for all $0\le d <d_0$ for some $d_0\ge 1$.
Let  \beqn
D(\hb ):=Z'_{d_0}-Z_{d_0}
&=& A\hb ^{-2r-1}+B\hb ^{-2r}+....\\
&=&\hb ^{-2r}(A/\hb  +B + O(\hb )) ,\eeqn
where $A, \ B \in H^*_{T\ti T'}(X)$. ($A$ might be $0$.)
This is possible since $<D,\phi _v>_0$ for all $v$ are polynomials
of $1/\hb$ over $H^*_{(T\ti t')}$
and so $D$ is a polynomial of $1/\hb$ over $H^*_{T\ti T'}(X)$.
Consider the coefficient of $q^{d_0}$ in $W(Z')-W(Z)$, which can
be set $\delta (D)
  =\int _{V}e^{(p+d_0\hb )z}D(\hb )+e^{pz}D(-\hb ) $.
If $r=0$, then $D=0$ since $D\equiv 0$ modulo $(1/\hb )^2$.
Assume $r\ge 1$. We shall show that $A=0=B$, which
implies by induction that $D=0$.
Notice that, since $\kd (D)$ is a polynomial of $\hb $,
\beqn
O(\hb ^2)=\hb ^{2r}\kd (D) &=&\int _{V}e^{(p+d_0\hb )z} (A/\hb  + B+ O(\hb )) 
         +e^{pz}(-A/\hb  +B +O(\hb )) \\
  &=&\int _{V}e^{pz}Ad_oz +2Be^{pz} +O(\hb ).
\eeqn
So,
\beqn
0 &=& d_0z\int _{V}e^{pz}A + 2 \int _{V}Be^{pz} \\
 &=& \sum _{v\in X^T}(d_0ze^{p_vz}A_v+ 2 e^{p_vz}B_v)\frac {1}{i_v^*(\phi _v)},
\eeqn
where $A_v$, $B_v$, and $Euler(V)_v$
are the restrictions of $A$, $B$, and $Euler(V)$ to the fixed
point $v$, respectively.
Since
$p_vz$ are different as $v$ are different (this can be seen in
section \ref{grass} and \ref{flag}),
$e^{p_v z}$ and $ze^{p_vz}$
are independent over  $H^*_{(T\ti T')}$.
So we conclude that
$A_v=0=B_v$ for all $v$, and hence $A=0=B$.

\subsection{Proof of theorem \ref{thmtr}}

It suffices  to show the transitivity of the action.
Let $Z_1$ and $Z_2$ be in class ${\mathcal P}({\mathcal C})$
and let $Z_1=1+o(1/\hb )$.

Since $\deg q\ge 0$, we may let
$$ Z_2=Z_2^{(0)} +Z_2^{(1)}\frac 1\hb  + o(\frac 1\hb ),$$ where
$Z_2^{(0)}\in H^*_{T\ti T'}(X)[[q]]$ 
is of  homogeneous of degree $0$ and $Z_2^{(1)}$
is homogeneous of degree $1$. Furthermore,
$Z_2^{(0)}(q)\in H^*_{T\ti T'}(X)[[q]]$ is a $q$-series with coefficients 
in $\QQ$, $Z_2^{(0)}(0)=1$, and
$Z_2^{(1)}(q)$ is a $q$-series with coefficients in $H^*_{T\ti T'} [p]$
by degree counting.

We may let $$\frac{Z_2^{(1)}(q )}{Z _2^{(0)}(q)}=
\sum _i (f_i(q)\cdot p_i) + g(q),$$ 
where $f_i(q)$ 
are  pure $q$-series over $\QQ$ of degree 0 and $g(q )$ are degree 1
in $H^*_{T\ti T'}[[q]]$.
In addition, $f_i(0)=0=g(0)$.

Now, consider operations on $Z_1$: first, coordinate changes,
\[ Z_1'=\exp (f(q)p/\hb )Z_1(q\exp (f(q)), \hb )
 = 1+ f(q)p/\hb  +o(1/\hb ), \]
second, multiplication by $\exp (g(q)/\hb )$,
\[Z_1''=\exp (g(q) /\hb  ) Z_1'
    = 1 + \frac 1\hb  (f(q)p + g(q)) + o(1/\hb ), \]
finally, multiplication by $Z_2^{(0)}(q)$,
\[ Z_1''' = Z_2 ^{(0)}(q)Z_1''
 =Z_2 ^{(0)} +\frac 1\hb  Z_2 ^{(1)} + o (1/\hb ).\]
According to the uniqueness lemma, the last one $Z_1'''$ must be equal to
$Z_2$ since
$Z_1'''\cong Z_2$ modulo $(1/\hb )^2$.

\subsection{Transformation from $J^V$ to $I^V$}
We explain the transformation introduced in the introduction. 
Let $\tilde{Z}$ be the nonequivariant specialization of $Z$.
Let $Z_1$ and $Z_2$ be in class ${\mathcal P}({\mathcal C})$
and let $Z_1=1+o(1/\hb )$.
Now let us specialize the equivariant setting to nonequivariant one.
Let $J^V=e^{(t_0+pt)/\hb }\tilde{Z_1}(q)$ and 
$I^V=e^{(t_0+pt)/\hb }\tilde{Z_2}(q)$.
Then, they are equivalent up to the unique coordinate
change $t_0\mapsto t_0 + f_0(q)\hb  + f_{-1}(q)$
and $t\mapsto t_i + f_i(q)$, $i=1,...,k$, where
$f_j\in \QQ [[q]]$ for all $j$,
$f_0$ and $f_i$ ($i=1,...,k$) have degree $0$, $f_{-1}$ has degree $1$; and
$f_j(0)=0$ for all $j$.

\section{The modified B series}
Let $X$ be a homogeneous manifold with
the torus $T\ti T'$ action. From now on
let $V=L_1\oplus ...\oplus L_l$ be an equivariant decomposable
convex vector bundle over $X$, where $L_i$ are 
line bundles.

\subsection{The correcting Euler classes}
Let $x=(x_1,...,x_l)$ be indeterminant.

Define a polynomial of $x$ over $\ZZ [\hb ]$ for $\beta\in \Lambda$:
$$ H_{\beta}(x,\hb )
:=\prod _{i=1}^{l}\prod _{m=0}^{<c_1(L_i),\beta >}(x_i+m\hb ).$$
Set $$H'_{\beta}(x,\hb ):=\frac{H_{\beta}(x,\hb )}{\prod x_i}. $$
We treat each linear factor $(x_i+m\hb)$ of $H_\beta$ as a Chern character.
Define \[ \Phi ^V (q,\hb ):=\sum _{d\in\Lambda }\sum _a
q^d<\frac{T_a}{\hb (\hb -c)}>^X_dT^a
E(H'_d(x,\hb ))(c_1(L),\hb ),\]
where $c_1(L)=(c_1(L_1),...,c_1(L_k))$.   

\bigskip

{\it Claim}
1. $(p_i)_w = (p_i)_v-<p_i, \beta _{v,w}>\ka _{v,w}$,

2. $c_1(L)_w = c_1(L)_v - < c_1(L), \beta _{v,w}>\ka _{v,w}$

3. $E(V_{v,w,m})=E(H_{m\beta _{v,w}})(c_1(L)_v,-\frac{\ka _{v,w}}{m}).$

{\it Proof:} Let $U$ be any equivariant convex line bundle.
On the ray $o(v,w)$ ($\cong \PP ^1$), 
we have a homogeneous coordinate $[z_0:z_1]$ such that the induced action
on the ray is linear (because of the equivariant embedding theorem).
We have also global sections $z_0^n,z_0^{n-1}z_1,...,z_1^n$ of the restriction
$U|_{\PP ^1}$ of $U$ to the ray,
where $[1:0]=w$, $[0:1]=v$ and $n=<c_1(U), \beta _{v,w}>$.
We know that $z_0^n$, $z_1^n$, $z_0/z_1$ have the characters $c_1(U)_w$,
$c_1(U)_v$ and $-\ka _{v,w}$, respectively. This concludes the proof.

\bigskip

The first one in the claim shows that we have the well-defined classes 
${\mathcal P}(X,V,E)$. 
(Otherwise, the mirror group transformation may not 
preserve the class ${\mathcal P}(X,V,E)$.)

\bigskip

{\theorem\label{thmmo}
Suppose $c_1(TX)-c_1(V)$ is in ample cone. Then
$\Phi ^V$ is in the class ${\mathcal P}(X,V,Euler)$.  }

\bigskip

Notice that for $\beta ' \le \beta$
 \begin{eqnarray}
H_\beta (x-<c_1(L), \beta ' >\hb , \hb )
   &= &H_{\beta '}(x,-h)H'_{\beta -\beta '}(x,\hb ) \label{Hdo} \\
   H'_\beta (x, \hb )
   &=&H'_{\beta '}(x,h)H'_{\beta -\beta '}(x+<c_1(L),\beta '>\hb ,\hb ) , 
\label{Hre}
\end{eqnarray}
which will show the polynomiality of double construction
and the almost recursion relation for $\Phi ^V$, respectively.

\subsection{The proof of theorem}

The homogeneous of $\Phi ^V$ is clear when $E=Euler$ 
and the rest properties will be proven for general $E$.

First of all, it is easy check to see 
\[ \Phi ^V \in H^*_{T\ti T'}[[\hb ^{-1}]][[q]] .\]
For the polynomiality, consider
\begin{eqnarray}
&& \int _V \Phi ^V(q,\hb ) e^{pz}\Phi ^V(qe^{-\hb z},-\hb ) \label{Hdoo} \\
&=& \sum _d \sum _{d^{(1)}+d^{(2)}=d, a} q^{d^{(1)}}
<\frac{T_a E(H_{d^{(1)}})(c_1(L),\hb )}{E(V)\hb (\hb -c)}>^X_{d^{(1)}} 
\nonumber \\
&&
 q^{d^{(2)}}<\frac{T^ae^{(p-d^{(2)})z} E(H_{d^{(2)}})(c_1(L),-\hb )}{
-\hb (-\hb -c)}>^X_{d^{(2)}}, \nonumber \end{eqnarray}
where $<T_a, T^b>_0^X=\delta _{a,b}$.
Let us use the notation and facts in \ref{pfdo}.
Since \beqn && E(H_{d^{(1)}})(c_1(L),\hb )E(H_{d^{(2)}})(c_1(L),-\hb ) \\ &=&
E(H_d(x-<c_1(L),d^{(2)}>\hb,\hb ))(c_1(L),\hb )E(V) \eeqn
 from (\ref{Hdo}),
the universal class $U(c_1(L))$ corresponding to $c_1(L)$ restricted to
$G_{d_1,d_2}(X)$ is
\[ c_1(L)-<c_1(L),d^{(2)}>\hb ,\]
and $e_1\circ \pi _1 = e_1\circ \pi _2$, 
(\ref{Hdoo}) is equal to 
\beqn
&&
\sum _d  q^d\int _{G_{d^{(1)},d^{(2)}}(X)}
      \frac{ e^{(\pi _2e^*_1p-d^{(2)}\hb)z}E(H_d)(U(c_1(L)),\hb )}
         {[N_{G_d(X)/G_{d^{(1)},d^{(2)}}(X)}]} \\    
&=& \sum _d q^d  \int _{G_d(X)} e^{Pz}E(H_d)(U(c_1(L)),\hb ), \eeqn
which shows the polynomiality.

Now let us check the almost recursion relation.
Let
$$S _v^X(q,\hb ):=
<S ^X,\phi ^X_v>^X_0=\sum _d S _{v,d}^X(\hb )q^d.$$ Since (if $d\ne 0$)
$$S _{v,d}^X(\hb )=\sum _{w\in o(v), 0<m; m\beta _{v,w}\le d}
\frac {C^X_{v,w,m}}{\hb (\ka _{v,w} + m\hb )}
S ^X_{w,d-m\beta _{v,w}}( -\frac{\ka _{v,w}}{m}) $$
and 
\[ E(H'_{\beta})(c_1(L)_v,-\frac{\ka _{v,w}}{m})
 =\frac{E(V_{v, w,m})}{E(V)_v}
E(H'_{\beta -m\beta _{v,w}})(c_1(L)_w, -\frac{\ka _{v,w}}{m}) \]
 from (\ref{Hre}) and the Claim,
we obtain that
\beqn
&&\Phi ^{V}_{v,d}(\hb ):=<\Phi ^V_d(\hb ),\phi _v>^V_0
= R_{v,d}  \\ &+&\sum _{w\in o(v), 0<m; m\beta _{v,w}\le d}
\frac{C^X_{v,w,m}E(V_{v,w,m})}
{E(V)_v\hb (\ka _{v,w} +m \hb )} \\
&& \ti \Phi ^X_{w,d-m\beta _{v,w}}
(-\frac{\ka _{v,w}}{m})
E(H'_{d-m\beta _{v,w}})(c_1(L)_w,-\frac{\ka _{v,w}}{m}),
\eeqn
where $R_{v,d}$ is indeed
 a polynomial of $1/\hb $ over $H^*_{(T\ti T')}$.

However, since
$$ C^V_{v,w,m} = \frac{C^X_{v,w,m}E(V_{v,w,m})}{E(V)_v},$$
$\Phi ^V_v(q,\hb )$ has the same almost recursion coefficients 
$C^V_{v,w,m}$ with $S^V$.

\subsection{Proof of main theorem \ref{thmmain}}
Recursion lemma \ref{lemmare} and double construction lemma \ref{lemmado}
show that $S^V$ is in class ${\mathcal P}(X,V,Euler)$. 
Certainly $S^V$ is form of
$1+o(\hb ^{-1})$. According to theorem \ref{thmmo}, $\Phi ^V$ also belongs to
${\mathcal P}(X,V,Euler)$. Then theorem \ref{thmtr} concludes the proof.
(We use the condition that $E=Euler$, in order to make sure that
$S^V$ and $\Phi ^V$ are homogeneous of degree 0.)

\section{Grassmannians}\label{grass}
\subsection{Notation}
Let $e_1,...,e_n$ form the standard basis of $\CC ^n$,
$T=(\CC ^\ti )^n$ the complex torus, and
$X:=Gr(k,n)$ the Grassmannian, the set of all
$k$-subspaces in $\CC ^n$.
As usual, let
 $T$ act on $Gr(k,n)$ by the diagonal action.
The fixed points $v=(i_1,...,i_k)$ are then the $k$-planes generated
by vectors $e_{i_1},...,e_{i_k}$.
Denote by $\CC ^n\ti X$  the trivial vector bundle with the
standard action.
Then we may consider
$L$, the determinant of the bundle dual to
the $T$-equivariant universal $k$-subbundle of $\CC ^n\ti X$.
Define
$V=L^{\ot l}$, $l>0$.
Denote
by $p$ the equivariant class $c_1(L)$.
We may identify $H^*(BT)$ with $\QQ [\ke _1,...,\ke _n]$
by the correspondence that $\ke _i$ is also
denoted the equivariant Chern class of the line
bundle over a point equipped with $T$ action as
the representation of the character $\ke _i$.
With respect to the Chern class of $L$,
we shall write $d\in \ZZ=H_2(X, \ZZ )$.
\subsection{A series}

\subsubsection{Fixed points}
Let $v$ be, say, $(1,2,..,k)$. Then around the point,
a local chart can be described by
\[\left( \begin{array}{cccc}
   1 & 0 & ... & 0\\
   0 & 1 &  & 0\\
   0 & 0 & ...& 1 \\
   x_{1,1} & x_{1,2} & ...& x_{n-k,k} \\
    & ...  &      \\
   x_{n-k,1} & x_{n-k,2} & ... & x_{n-k,k}
\end{array}\right).\]
For each complex value $(x_{i,j})$ the column vectors
in the matrix span a $k$-plane which stands for
a point in $Gr(k,n)$.
Then in the chart the action by $(t_1,...,t_n)\in T$ is described as
follows:
\[\left(\begin{array}{cccc}
 x_{1,1} & x_{1,2} & ...& x_{n-k,k} \\
    & ...  &      \\
   x_{n-k,1} & x_{n-k,2} & ... & x_{n-k,k}
\end{array}\right)
\mapsto \left(\begin{array}{cccc}
t_1^{-1}t_{k+1} x_{1,1} & t_2 ^{-1}t_{k+1} x_{1,2}
& ...& t_k^{-1}t_{k+1} x_{1,k}\\
 & & ...& \\
t_1^{-1}t_n x_{n-k,1} & t_2 ^{-1}t_{n} x_{n-k,2} & ...
& t_k^{-1}t_{n}x_{n-k,k}
\end{array}\right).\]

In each isolated fixed point of the Grassmannian there is
$\dim Gr(k,n)$-many 1-dimensional orbit (ray) passing through
the point. For instance, if $v=(1,2,...,k)$, then
there is only one ray $(v,w)$
 from $v$ to $w=(...,\hat{i},..., j)$
for any $i\le k <j$. These rays have
degree $1\in H_2(X, \ZZ )$.

\subsubsection{The Euler classes}
Notice that the tangent space at $v=(1,2,...,k)$ of the ray
connecting $v$ to $w=(...,\hat {i},...,j)$
has the character $\ka (v,w)=\ke _j - \ke _i$, where $j>k\ge i$.
Similarly one can find out the characters for the other
cases.

Let $f:\PP ^1\ra Gr(k,n)$ be a $m$-fold morphism totally ramifying
the ray over  $v$ and $w$.
The $T$ representation space $H^0(f^*L^{\ot l})$ has
the orbi-characters 
\beqn 
\frac {ap_v + b p_w}{m}
= lp_v -\frac{\ke _j-\ke _i}{m}b, 
 \text{ for } a+b=lm, \ a\ge 0,\  b \ge 0, 
\eeqn
where $p_v=-(\ke _1+...+\ke _k)$ 
and $p_w=-(\ke _1+...+\hat{\ke _i} +....+\ke _k + \ke _j)$ 
are $p=c_1(L)$ restricted to the fixed
points $v$ and $w$, respectively.

\section{The flag manifolds}\label{flag}
We analyze
fixed points of the maximal torus actions and the invariant
curves connecting two fixed points. This explicit description would be
useful also to find $S^X$ explicitly.
\subsection{The complete flag manifolds}
Let $X$ be the set of all Borel subgroups of a simply connected
semi-simple Lie group $G$. It is a homogeneous space with
the $G$-action by conjugation.
Then the maximal torus $T$-action (---fix one---)
has isolated fixed points. They are exactly Borel subgroups
containing $T$. The fixed points are
naturally one-to-one corresponding to the set of 
Weyl chambers. Each Borel subgroup containing $T$ gives rise to
a negative roots (---our convention---)
of $B$ and so a chamber associated to the positive
roots. Let $C$ be the set of chambers.
The tangent line subspace associated to
the positive root $\ka$ has the character $\ka$.

There is, if one fix a fixed point $v$, 
a natural correspondence between the
$H^2(X,\ZZ )$ and the characters of $T$.
Then the K\"ahler cone is exactly the 
positive Weyl chamber $v$. Notice that
the fundamental roots span the  K\"ahler cone.
Consider co-roots $\ka ^\vee$. They span
the Mori cone. We can identify the Mori cone $\Lambda$
with the non-negative integer span of
co-roots.
\subsection{The generalized flag manifolds}
Let $X$ be the set of all parabolic subgroups with a
given conjugate type. Let $T$ be a maximal torus of $G$.
Then the fixed loci are isolated fixed points
consisting of  parabolic subgroups containing $T$.

\subsubsection{Rays}\label{ray}
Let us choose a fixed point $P\supset T$. Then the rays
at the fixed points are described by the following way.
(The rays are by definition the 1-dimensional orbits of $T$
passing through $P$.)
Fix $B$ a Borel subgroup in $P$ containing $T$.
First consider the $T$-equivariant fibration,
$G/B \ra G/P$ and
the rational map 
to $G/B$ by $\exp (zX_\ka )\in G$, $z\in \CC$,
where $X_\ka$ is an eigenvector of the positive root
$\ka$.
Since $\exp H\exp (zX_\ka )\exp -H
=\exp (z\exp\ad H (X_\ka ))=\exp (z\exp (\ka (H)) X_\ka )$ for
$H\in \Lie T$, 
we conclude that it is a
$T$ - invariant stable map.
By the composition of the fibration, we obtain all
the rays. They are effectively labeled by the positive roots which
are not roots of $P$. So there are exactly $\dim X$-many
rays at each fixed point. The tangent line at the ray
has the character $\ka$.

\subsubsection{The K\"ahler cone}
Here we need the Levi-decomposition of $P$ and then
consider simple roots $\{
\ka _i\} _{i\in P(\Delta )}$ which
are not roots of the semi-simple part of $P$.
Then the fundamental roots with respect to $P$ is, by definition,
$\{ \kl _i \} _{i\in P(\Delta )}$, where $\kl _i$ are dual to
$\ka _i^\vee $.

Choose a fixed point $P$.
We may identify  $H^2(X,\ZZ )$ with the set of 
integral weights according to Borel-Weil theorem.
Then the K\"ahler cone is the set of all 
dominant integral weights {\em with respect to} $P$.

\subsubsection{Homogeneous line bundles}
One can produces all very ample line bundles by  homogeneous
line bundles associated to irreducible representations of $P$
with  highest weights $\kl$. The weights corresponding to
the very ample line bundles are exactly the positive integral combination
of the fundamental weights with respect to $P$.
We shall denote by ${\mathcal O}(\kl )$
the homogeneous line bundle associated to
the (1-dimensional) highest weight  $\kl =\sum _{i\in P(\Delta )}
a_i\kl _i$ representation of $P$. 
It is a  very ample bundle if and only if
$a_i >0$ for all $i$.

This also shows that the ray $\PP ^1$ associated to $\ka$ has
the homology class \lq\lq $\ka ^\vee$," in the sense that 
$<\PP ^1,c_1({\mathcal O}(\kl )>=(\ka ^\vee ,\kl )$.
{\em We shall use $\ka ^\vee$ to denote the homology class.}

\subsubsection{$\sum (p_i)_vz_i$ are different for different $v$.}
Consider a line bundle $L$ associated
to $\kl =\sum _{i\in P(\Delta )}a_i\kl _i$. (Here in advance,
we have to fix $P\supset T$.) 
Let $S_{\ka}$ denote the Weyl group element of the
reflection associated to the root $\ka$. 
Then the line bundle
is $L={\mathcal O}(S_{\ka}(\kl ))$ if one look at it
with respect to another \lq\lq origin " 
$P'=\exp (\frac{\pi}{2}(X_\ka -Y_\ka)) 
P \exp (-\frac{\pi}{2}(X_\ka -Y_\ka))$ where
$[X_\ka ,Y_\ka]=H_\ka$, $[H_\ka ,X_\ka ]=2X_\ka$ and
$[H_\ka ,Y_\ka]=-2Y_\ka$.
This $P'$ is the other $T$-fixed point lies in the ray
associated to $\ka$ which passes through $P$.
(Because of the $SL(2,\CC )$-equivariant map from
$\PP ^1$ to the ray, it is enough to check it
when $G=SL(2,\CC )$, which is obvious.)

\subsubsection{$V_{v,w,m}$ and $N_{v,w,m}$}
Let $V={\mathcal O}(\kl )$.
Let $\psi :\PP ^1\ra X$ be a stable map totally ramifying one of
rays, passing through $P\supset T$. Suppose the
ray is associated to a positive root $\ka $ with respect to $P$
and $f$ is a $m$-multiple branched
cover representing an isolated
$T$-fixed point of $\overline{M}_{0,0}(X,m\ka ^\vee)$. 
Then the $T$-representation
space $H^0(\PP ^1, f^*({\mathcal O}(\kl ))$ has the characters,
\[ \kl - a\frac{\ka}{m}, \ a=0,...,m(\kl ,\ka ^\vee). \]
To see it, use the coordinate $z\in \CC$ around the fixed
point and $\exp H\exp (zE_\ka )\exp -H
=\exp (z\exp\ad H (E_\ka ))=\exp (z\exp (\ka (H)) E_\ka )$.
Similarly,
$N_{v,w,m}$ has characters
\beqn
 \delta -a\frac{\ka }{m} && \text{ for } \ka \ne \delta >0,
a=0,...,m(\delta ,\ka ^\vee ),  \\
 \ka - a\frac{\ka }{m} && \text{ for }  a=0,...,\hat{m},...,2m,
\eeqn
where $\delta >0$ means $\delta$ is a positive root with respect to $P$.

{\it e-mail address}: {\tt bumsig@math.ucdavis.edu }


\begin{thebibliography}{9}

\bibitem{BCKV1} V. Batyrev, I. Ciocan-Fontaine,
               B, Kim and D. van Straten, {\it
Conifold transitions and mirror symmetry for complete
intersections in Grassmannians, } Preprint 1997.

\bibitem{BCKV2}
          --------, {\it Mirror symmetry and toric degenerations
of partial flag manifolds,} In preparation.

\bibitem{BV} V. Batyrev and D. van Straten,
{\em Generalized hypergeometric functions
and rational curves on Calabi-Yau complete intersections in toric
varieties}, Commun. Math.  Phys. {\bf 168} (1995), 493-533,
(alg-geom/9307010).

\bibitem{Be} K. Behrend, 
 {\it Gromov-Witten invariants in algebraic geometry},
      To appear in Inventiones Math..
\bibitem{BF} K. Behrend and B. Fantechi, 
 {\it The intrinsic normal cone,} To appear in Inventiones Math..
 
\bibitem{BM} K. Behrend and Yu. Manin,
   {\it Stacks of stable maps and Gromov-Witten invariants},
    Duke Math. J. vol 85, No1 (1996) 1-60.



\bibitem{Ca} P. Candelas, X. de la Ossa, P. Green
  and L. Parkes,
{\it A pair of Calabi-Yau manifolds as an exactly soluble 
superconformal theory,} Nuclear Phys. B {\bf 359} (1991), 21-74.

\bibitem{Ionut} I. Ciocan-Fontaine, {\it
   On quantum cohomology rings of partial flag varieties,}
   Institut Mittag-Leffler Report No. 12, 1996/1997.


\bibitem{Du} B. Dubrovin, {\it The geometry of 2D topological field
theories} in Integrable systems and quantum groups, Lecture Notes in
Math. {\bf 1620}, Springer-Verlag, Berlin, 1996, 120-348.
 


\bibitem{FuP} W. Fulton and R. Pandharipande, 
{\it Notes on stable maps and quantum cohomology,}
 To appear in  Proc. AMS summer school, 1995,
Santa Cruz.

\bibitem{GH} A. Givental, {\it
      Homological geometry I: Projective hypersurfaces,}
     Selecta Math., New Series, {\bf 1} (1995), No. 2, 325-345.

\bibitem{GE} --------, {\it Equivariant Gromov-Witten invariants,}
       IMRN 1996, No. 13, 613-663.
 
\bibitem{GS} --------, {\it Stationary phase integrals,
quantum Toda lattices, 
flag manifolds and the mirror conjecture } Preprint 1996.

\bibitem{GT} --------, {\it A mirror theorem for Toric complete
             intersections}, Preprint 1997.
 


\bibitem{KT}  B. Kim,
{\it Quantum cohomology of flag manifolds $G/B$ and
quantum Toda lattices}, Preprint 1996.
\bibitem{Ko}  M. Kontsevich, 
{\it Enumeration of rational curves via torus
actions}, in  The moduli space of curves, 
R. Dijkgraaf, C. Faber, and
G. van der Geer, (Eds.),
Progr. Math. 129, Birkh\"auser, Boston, 1995, 335-368.

\bibitem{LT}  J. Li and G. Tian, {\it Virtual moduli cycle and
Gromov-Witten invariants of algebraic varieties,} Preprint 1996.



\bibitem{Th} J. Thomsen, {Irreducibility of
       $\overline{M}_{0,n}(G/P,\beta )$,}
Institut Mittag-Leffler Report No. 8, 1996/97.
 
\bibitem{To} E. Tj\o tta, {\it Rational curves on the
space of determinantal nets of conics,} Thesis 1997, 
University of Bergen.

\bibitem{W} E. Witten, {\it Two-dimensional gravity and
 intersection theory on moduli space,}
 Survey in Differential Geometry {\bf 1} (1991) 243-410.

\end{thebibliography}
\end{document}